\documentclass[10pt]{article}

\twocolumn

\usepackage{amssymb,amsfonts,amsmath,epsfig,graphicx,float, wrapfig}
\newcommand{\sech}{\text{sech} \, }
\newcommand{\arsinh}{\text{arsinh} \, }

\newcommand{\R}{\mathbb{R}}
\newcommand{\C}{\mathbb{C}}
\newcommand{\Z}{\mathbb{Z}}

\newcommand{\cS}{\mathcal{S}}

\newtheorem{example}{Example}

\begin{document}

\author{Hans Engler\\
Department of Mathematics and Statistics\\
Georgetown University\\
Washington, D.C. 20057, U.S.A.
}

%
%
%


\date{\today}
\title{Computation of Scattering Kernels in Radiative Transfer}

\maketitle

\section*{Abstract}
\small{
This note proposes rapidly convergent computational formulae for evaluating scattering kernels from radiative transfer theory. The approach used here does not rely on Legendre expansions, but rather uses exponentially convergent numerical integration rules. A closed form for the \textsc{Henyey-Greenstein} 
scattering kernel in terms of complete elliptic integrals is also derived. 
}

\medskip
\noindent
\small{\textbf{Keywords.} Plane-parallel scattering, scattering kernel, phase function, trapezoidal rule, complete elliptic integral.}

\section{Introduction and Background}

In theories of radiative transfer and of neutron transport, the interaction between radiation and a scattering medium is described by a \emph{phase function} $p:[-1,1] \to \R^+$ with the property $\int_{-1}^1 p(t) dt = 1$. Let $S^2$ denote the unit sphere in $\R^3$. In a single scattering event, a photon or neutron that arrives from a given direction  $\Theta \in S^2$ is scattered into the direction 
$\Theta' \in S^2$ according to the probability density $\Theta' \mapsto \frac{1}{4\pi} p(\Theta \cdot \Theta')$. This probability density then appears as an integral kernel in an integro-differential equation that describes multiple scattering. In the special case of plane-parallel scattering, one rewrites this equation in spherical coordinates and expands it as a Fourier series in the azimuthal angle $\phi$.  
This requires the evaluation of 
\begin{equation}
c_m \int_0^{2 \pi}  p(\cos \theta \cos \mu + \sin \theta \sin \mu\cdot \cos \phi') \cos m \phi' \, d \phi'
\label{eq_average}
\end{equation}
for $0 \le \theta, \, \mu \le \pi, \, m = 0, \, 1, \, 2, \dots$, where $c_m = \frac{1}{2\pi}$ for $m > 0$ and $c_0 = \frac{1}{\pi}$.
Using the notation $x = \cos \theta, \, y = \cos \mu$, one therefore arrives at the problem of evaluating the \emph{scattering kernels} $P_m(x,y)$, defined for $-1 \le x, \, y \le 1$ and for $m = 0, \, 1, \dots$ as
\begin{equation}
c_m \int_0^{2 \pi}  
p(xy  + \sqrt{1-x^2}\sqrt{1-y^2}\cos s) \cos m s \, d s
\label{eq_matrix}
\end{equation}
Let $L_n$ denote the n-th Legendre polynomial. From spherical harmonics one obtains that the Legendre expansion of the phase function $p$ directly 
leads to expansions for the $P_m$. If $p(x) = \sum_{n = 0}^\infty \alpha_n L_n(x)$ for all $x \in [-1,1]$, then in particular
\begin{equation}
P_0(x,y) = \sum_{n = 0}^\infty \alpha_n L_n(x)L_n(y) 
\label{eq_LegendreExp}
\end{equation}
for all $-1 \le x, y \le 1$. The higher order terms $P_m, \, m > 0$ can be expanded into associated Legendre functions, using again the Legendre coefficients of $p$. The classical reference \cite{chandrasekhar1960} contains a mathematical presentation of the theory of radiative transfer. A modern account with more physical details may be found in \cite{liou2002}. 

While the evaluation of eq.~(\ref{eq_LegendreExp}) and its higher order versions  is in theory straightforward, there are several practical difficulties. Firstly, the formula requires knowledge of the Legendre expansion of $p$. While the methods to find $p$ for a particular scattering medium indeed produce Legendre expansions (see \cite{wiscombe1980}), it may be desirable to have more compact representations of a scattering function, in which case the Legendre expansion is not readily available and not easy to compute (there is no ``fast Legendre transform"). Secondly, for cases of strongly forward-peaked scattering, several hundred terms of the Legendre expansion may be needed to evaluate each $P_m$ even to modest accuracy. This requires care when evaluating the Legendre polynomials in eq.~(\ref{eq_LegendreExp}), and it is computationally expensive in any case.    

In this note, the direct numerical integration scheme known as trapezoidal rule is proposed to evaluate $P_m$ from eq.~(\ref{eq_matrix}), as an alternative to the Legendre expansion in eq.~(\ref{eq_LegendreExp}). The method is known (\cite{Trefethen2014}) to converge exponentially if the integrand is analytic. This highly desirable property is exploited systematically in this note.
The relation between the convergence rate and the location of the singularities (points or regions of non-analyticity) of the phase function is explained. The second contribution of this note is the derivation of a closed form of the scattering kernel $P_0$, in the special case of the Henyey-Greenstein phase function (\cite{henyey1941}). It expresses $P_0$ in terms of a complete elliptic integral and can be evaluated very rapidly without any expansions or numerical integration, even for cases of extremely forward-peaked scattering. Numerical examples are given to demonstrate the approach.

\subsection{Henyey-Greenstein Scattering Function}
\begin{figure}[ht]
\begin{center}
\resizebox{3in}{!}{\includegraphics{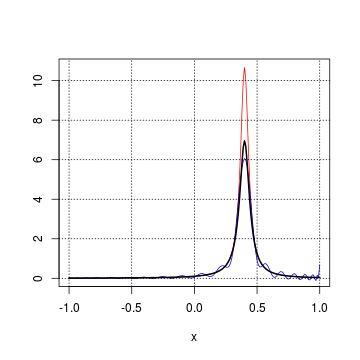}}
\caption{Henyey-Greenstein scattering kernel $H(x,y_0;g)$  for $g = 0.95,  \, y_0 = 0.4, \,  -1 \le x \le 1$. Black: exact evaluation using eq.~(\ref{eq_HG_exact}). Blue: Eq.~(\ref{HG0}) with $N = 40$ terms. Red: Eq.~(\ref{eq_trapezoidal}) with $N = 40$ terms.}
\label{fig-HG1a}
\end{center}
\end{figure}

The Henyey-Greenstein phase function (\cite{henyey1941}) \footnote{It appears that Chandrasekhar was not aware of this work when he wrote his classical treatise 
\cite{chandrasekhar1960} in 1950.} was proposed to describe interstellar scattering and is given by the formula
\begin{equation}
p_{HG}(x,g) = \frac{1}{2} \frac{1-g^2}{(1 + g^2 - 2g x)^{3/2}}
= \sum_{\ell = 0}^\infty \frac{2 \ell + 1}{2}g^\ell L_\ell(x)
\label{HG}
\end{equation}
where $-1 < g < 1$ is known as the asymmetry factor. This phase function has since been used in areas as diverse as scattering in cloudy and hazy atmospheres (\cite{hansen1969}), light scattering in seawater (\cite{haltrin2002}) and in tissue (\cite{niemz2007}), and even in computer graphics.  Then  eq.~(\ref{eq_matrix}) together with eq.~(\ref{eq_LegendreExp}) lead to the problem of evaluating  
\begin{equation}
H(x,y;g) =  \sum_{\ell = 0}^\infty \frac{2 \ell + 1}{2}g^\ell L_\ell(x) L_\ell(y)
\label{HG0}
\end{equation}
for $-1 \le x,y \le 1$. The series may be evaluated by using the first $N$ terms. 
Since $|L_\ell(x) L_\ell(y)| \sim \frac{2}{2 \ell + 1}$ with indeterminate sign, the direct evaluation of the series leads to problems when $g$ is close to 1, because then the series converges very slowly. This is illustrated in fig.~\ref{fig-HG1a}. 

\begin{figure}[ht]
\begin{center}
\resizebox{3in}{!}{\includegraphics{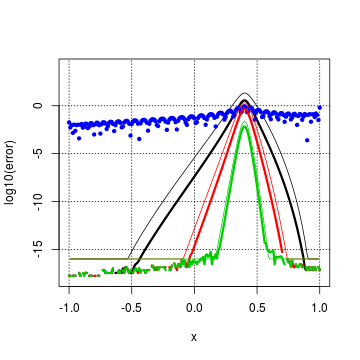}}
\caption{Log Errors for evaluating Henyey-Greenstein scattering kernel in fig.~\ref{fig-HG1a}. Blue: Eq.~(\ref{HG0}) with $N = 40$ terms. Black, red, green:  Eq.~(\ref{eq_trapezoidal}) with $N = 40, \, 80, \, 160$ terms. Thick: Computed errors. Thin: Error estimates from eq.~(\ref{eq_traperror}).}
\label{fig-HG1b}
\end{center}
\end{figure}

\section{An Exact Formula}  

We start with the product formula (see 18.17.6 in \cite{nist2010}) for Legendre polynomials
\begin{eqnarray*}
L_n(\cos \theta) L_n(\cos \mu) &=& \frac{1}{\pi} \int_0^\pi L_n (\cos \theta \cos \mu \\
& \; & + \sin \theta \sin \mu \cos s ) ds \, .
\end{eqnarray*}
Let 
\begin{equation}
H_0(x,y;g) =  \sum_{\ell = 0}^\infty g^\ell L_\ell(x) L_\ell(y) \, .
\label{HG00}
\end{equation}
Using the generating function for Legendre polynomials (see 18.12.11 in \cite{nist2010}), we therefore obtain
\begin{eqnarray}
& \,& H_0(\cos \theta, \cos \mu; g) \\
&=& \frac{1}{\pi} \int_0^\pi \sum_{k = 0}^\infty g^n L_n\left(\cos \theta \cos \mu + \sin \theta \sin \mu \cos s \right) \, ds\\
&=& \frac{1}{\pi} \int_0^\pi \frac{ds}{\sqrt{1 - 2g\left(\cos \theta \cos \mu + \sin \theta \sin \mu \cos s \right) + g^2}}
\, .
\end{eqnarray}
By eq.~(\ref{eq_K1}) with 
\begin{eqnarray*}
\alpha &=& 1 + g^2 - 2g\cos \theta \cos \mu \\
\beta &=& 2 g \sin \theta \sin \mu \\
 \alpha + \beta &=& 1 + g^2 - 2g\cos(\theta + \mu)
\end{eqnarray*}
this implies
\begin{equation}
H_0(\cos \theta, \cos \mu;g) =
\frac{2}{\pi\sqrt{\alpha + \beta}} K_0\left(\frac{2 \beta}{\alpha + \beta} \right) 
\label{eq_HSum1} 
\end{equation}
where $K_0$ is the complete elliptic integral of the first kind defined in eq.~(\ref{eq_K0}). 
Setting 
\begin{equation}
u_\pm = 1 + g^2 - 2g\cos(\theta \pm \mu)
\label{eq-upm}
\end{equation} 
this may also be written as
\[
H_0(\cos \theta, \cos \mu;g)
= \frac{2}{\pi\sqrt{u_+}} K_0\left(\frac{u_+ - u_-}{u_+} \right) \, .
\] 
This formula is closely related to formula  (5.10.2.1)\footnote{This was pointed out to me by an anonymous contributor on \texttt{math.stackexchange.com}.} in \cite{prudnikov1986} which has the form
\begin{equation}
H_0(\cos \theta, \cos \mu;g) = \frac{4}{\pi(\sqrt{u_+} + \sqrt{u_-})} K\left(\frac{\sqrt{u_+} - \sqrt{u_-}}{\sqrt{u_+} + \sqrt{u_-}}\right)\, .
\label{eq_HSum2} 
\end{equation} 
From eq.~(\ref{eq_HSum1}) and eq.~(\ref{eq-dK}) we obtain finally 
\begin{eqnarray}
&\,& H(\cos \theta, \cos \mu ; g)\\
 &=& \left(g \frac{d}{dg} + \frac{1}{2} \right) 
\frac{2}{\pi\sqrt{u_+}} K_0\left(\frac{u_+ - u_-}{u_+} \right)\\
&=& \frac{(1-g^2)}{\pi 
u_- \sqrt{u_+}} E_0\left(\frac{u_+ - u_-}{u_+} \right)   
\label{eq_HG_exact}
\end{eqnarray}
where $E_0$ is the complete elliptic integral of the second kind defined in eq.~(\ref{eq_E0}). In terms of $x, \, y$, this becomes
\begin{equation}
H(x,y;g) = \frac{(1-g^2)}{\pi 
w_- \sqrt{w_+}} E_0\left(\frac{4g\sqrt{1-x^2}\sqrt{1-y^2}}{w_+} \right)
\label{eq_HG_exactxy}
\end{equation}
with $w_\pm = 1 + g^2 - 2g \left(xy \mp \sqrt{1-x^2}\sqrt{1-y^2}\right)$. The formula may be used to evaluate $H(\cdot, \cdot; g)$ reliably even if $g$ is extremely close to 1. If $g = 1 - \varepsilon$ and $\varepsilon_0 \approx 10^{-16}$ is ``machine epsilon" in IEEE arithmetic, then $H(x,y;g)$ can be evaluated to relative accuracy $\varepsilon_0/\varepsilon$. 

\section{Fast Evaluation}

We now turn to the general case. For a given phase function $p$ and given $x, \, y \in [-1,1], \, m \in \{0, \, 1, 2, \dots\} $, we need to evaluate the integral given by eq.~(\ref{eq_matrix}).  Let 
\begin{eqnarray}
h_m(z) &=& p(A + B \cos z)\cos mz \\
 A &=& xy, \, B = \sqrt{1 - x^2}\sqrt{1 - y^2} \, .
\label{eq_def_h}
\end{eqnarray}
Note that $|A \pm B| \le 1$, with equality in one of the cases where $x = \pm y$.
The function $h_m$  is $2 \pi$-periodic on $\R$. It is known (\cite{Trefethen2014}) that if a function $f$ is periodic and analytic in a strip about the real axis, then the trapezoidal rule for approximating the integral $\int_0^{2\pi} f(t) dt$ converges exponentially fast. More precisely, assume that $f$ is $2\pi$-periodic and analytic in the strip $\cS_\alpha = \{z \in \C \, | \, |\Im z | < \alpha\}$ and satisfies $|f(z)| \le M$ for some constant $M \ge 0$ there. 
Choose a positive integer $N$ and set
\begin{equation}
I_N = \frac{2 \pi}{N} \sum_{k = 1}^N f\left( \frac{2 \pi k}{N} \right)\, , 
\label{eq_trapezoidal}
\end{equation}
then
\begin{equation}
\big| \int_0^{2 \pi} f(t) dt - I_N \big| 
\le \frac{4 \pi M}{e^{\alpha N} -1} = e^{-\alpha N} \frac{4 \pi M}{1 - e^{-\alpha N}} \, . 
\label{eq_traperror}
\end{equation}
The sum is just an approximation of the integral with the composite trapezoidal rule. The rate of convergence thus is much faster than for ordinary smooth (twice differentiable) functions, for which the error estimate has the form 
\begin{equation}
\big| \int_0^{2 \pi} f(t) dt - I_N \big| 
= \frac{2 \pi^3 |f´´(\xi)| }{3 N^2} 
\label{eq_traperrorC2}
\end{equation}
for some $\xi \in (0, 2 \pi)$. 

To use this result in the computation of eq.~(\ref{eq_matrix}), note first that for any $\alpha \in \R$, the function that takes $u + i \cdot v = z$ to 
$$
A + B \cos z =  A + B  
\cos u \cosh v + B \sin u \sinh v \cdot i 
$$ 
maps the strip $\cS_\alpha$ to an ellipse about $[A-B,A+B]$ in the complex plane which has focal points $A \pm B$ and major and minor axes with lengths $B \cosh \alpha$ and $B \sinh \alpha$.  
Therefore, if the phase function $p$ is analytic in a neighborhood surrounding the set  $[-1,1] \subset \C$, 
then $h_m$ is analytic in a suitable strip $\cS_\alpha$ with $\alpha > 0$. The domain of analyticity of $p$ is automatically symmetric with respect to the real axis. It should be emphasized that it is of course not necessary to determine this domain in order to use eq.~(\ref{eq_trapezoidal}).   

\begin{figure}[ht]
\begin{center}
\resizebox{3in}{!}{\includegraphics{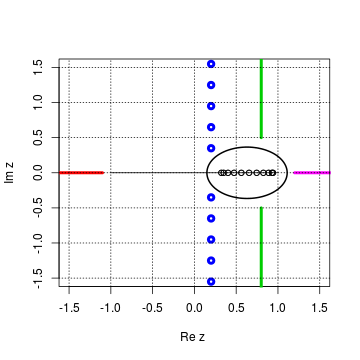}}
\caption{Integration domain for a phase function. }
\label{fig-HG2a}
\end{center}
\end{figure}

For an illustration, refer to fig.~\ref{fig-HG2a}. The plot shows the domain of an assumed phase function that is originally defined on the interval $[-1,1]$ (thin horizontal black line) and that can be extended into the complex plane (everywhere except at singularities shown as colored lines and circles). For a particular choice of $x,\, y$, the integral in eq.~(\ref{eq_matrix}) extends over $[A -  B, A+B] \subset [-1,1]$ (black circles). A suitable strip $\cS_\alpha$ is mapped to the ellipse surrounding this set where the phase function is analytic. Consequently, the trapezoidal approximation converges exponentially with a rate given by eq.~(\ref{eq_traperror}). The rate of convergence depends on $\alpha$ which in turn comes from the location of $[A-B,A+B]$ relative to the set of singularities of the phase function.    

\begin{figure}[ht]
\begin{center}
\resizebox{3in}{!}{\includegraphics{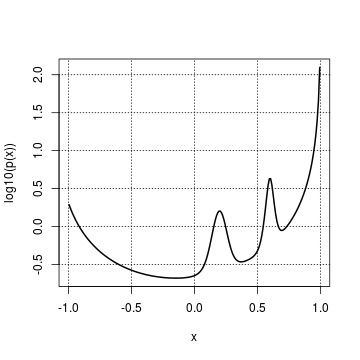}}
\caption{Multimodal phase function given by
eq.~(\ref{eq_multimodal}).}
\label{fig-HG2b}
\end{center}
\end{figure}

\begin{example} \textbf{Henyey-Greenstein Phase Function}. Consider a phase function that is analytic in $\C$ minus the ray $[\frac{1 + g^2}{2g}, \infty)$. Then the integrand in eq.~(\ref{eq_matrix}) is analytic in the strip defined by 
$$
\cosh \Im z < \left| \frac{(1+g^2)/2g - A}{B} 
\right| =  \left| \frac{(1+g^2)/2g - xy}{\sqrt{1-x^2}\sqrt{1-y^2}} 
\right| \, .
$$    
\end{example}
This case is illustrated by the magenta line in fig.~\ref{fig-HG2a}.  Similar comments apply to the case where $g < 0$ (red line).  The Henyey-Greenstein phase function defined in eq.~(\ref{HG}) is of this form. The error behavior of the trapezoidal rule~(\ref{eq_trapezoidal}) is illustrated in fig.~\ref{fig-HG1b} which shows the log-errors (thick lines) when this rule  is used to approximate the integral in eq.~(\ref{eq_matrix}) for $x \in [-1,1]$ and for fixed $y = y_0, \, g$, for three different choices of $N$. The plot shows that errors are maximal for $x \approx y_0$ and that the errors have the slowest decrease near this value as $N$ increases. Also plotted (thin lines, same colors)  are error estimates from eq.~(\ref{eq_traperror}), using $M = 1$ for simplicity. It is seen that the actual errors track the estimate very closely.     

\begin{example} Consider a phase function that is analytic except possibly on rays 
$x_0 \pm i  \, \delta$ to $x_0 \pm i  \, \infty$, where $x_0 \in \R, \, \delta > 0$. 
Then the integrand in eq.~(\ref{eq_matrix}) is analytic in the strip defined by 
\begin{equation}
- S < \sinh \Im z < S
\end{equation}
where  $S$ is the unique positive solution of the equation
$$
\frac{(A - x_0)^2}{S^2 + 1} + \frac{\delta^2}{S^2} = B^2 \, .
$$
\label{ex_bump}
\end{example}
To see this, find $\alpha$ such that the ellipse parametrized by 
$$
s \mapsto A + B (\cos s \cosh \alpha + i \cdot \sin s \sinh \alpha)
$$ 
passes through the points $x_0 \pm i \cdot \delta$, and set $S = \sinh \alpha$.  

To obtain examples, fix $\delta > 0, \, x_0 \in \R, \gamma > 0$ and consider functions of $x \in \R$
\begin{eqnarray}
f_1(x;x_0,m,\gamma) &\propto&  \left(1 + \left(\frac{x-x_0}{\delta} \right)^2 \right)^{- \gamma} 
\label{eq_f1}
\\
f_2(x;x_0,m) &\propto& \sech \left(\frac{x-x_0}{\delta} \right) 
\label{eq_f2}
\end{eqnarray}   
where the proportionality constants are chosen such that the integrals over $[-1,1]$ equal 1. Both functions have single maxima (peaks) at $x = x_0$ and the width of the peak is proportional to $\delta$. The Legendre expansions of these functions are generally not available in closed form. The function $f_1$ is analytic in the complex plane minus branch cuts from $x_0 \pm i  \, \delta$ to $x_0 \pm i  \, \infty$ (green lines in fig.~\ref{fig-HG2a}). The function $f_2$ is analytic in the complex plane minus poles at $z = x_0 \pm i \, \delta\left(\pi/2 + n  \pi \right), \; n \in \Z^+$ (blue circles in fig.~\ref{fig-HG2a}).  Therefore if $p = f_1$ or $p = f_2$, then the integrand in eq.~(\ref{eq_matrix}) is analytic in any strip $\cS_\alpha$ where $\alpha = \arsinh S$ and $S$ is as in example~\ref{ex_bump}.

The reader may note that the integrand in eq.~(\ref{eq_matrix}) contains factors $\cos m z$. These terms grow like $e^{ m |\Im z|}$ away from the real axis and their second derivatives contain the factor $m^2$. When the integrand is analytic and eq.~(\ref{eq_traperror}) can be used, the error is therefore proportional to 
$e^{-\alpha N} e^{\alpha m} M = e^{-\alpha (N-m)}M$. Thus the additional factor $\cos m z$ has the same effect as using $N-m$ points instead of $N$ points for the evaluation of eq.~(\ref{eq_trapezoidal}), resulting in a modest loss of accuracy. On the other hand, if the integrand is merely twice differentiable, the error from eq.~(\ref{eq_traperrorC2}) becomes proportional to $m^2/N^2$. Thus the additional factor $\cos m z$ now has the same effect as using only $N/m$ points instead of $N$ points, leading to a much larger loss of accuracy. This illustrates the powerful effect of having an analytic integrand. 

\subsection{Multimodal Phase Functions} 
\begin{figure}[ht]
\begin{center}
\resizebox{3in}{!}{\includegraphics{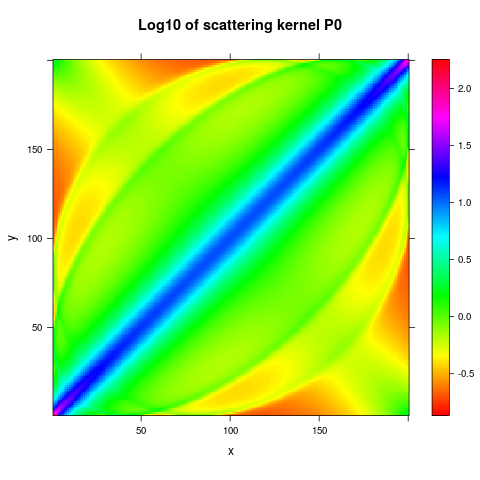}}
\caption{Logarithmic heat map of scattering kernel $P_0$ for the phase function given by eq.~(\ref{eq_multimodal}).}
\label{fig-HG3}
\end{center}
\end{figure}

\begin{figure}[ht]
\begin{center}
\resizebox{3in}{!}{\includegraphics{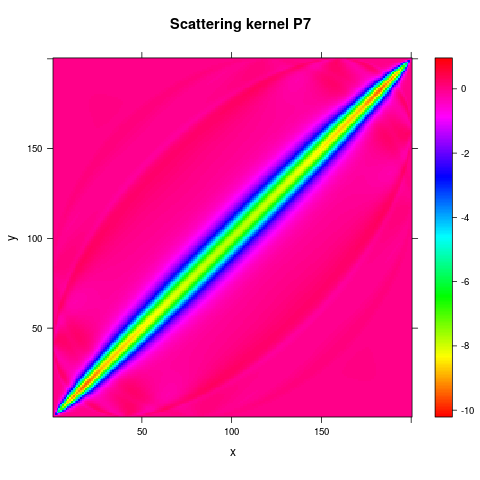}}
\caption{Heat map of scattering kernel $P_7$ for the phase function given by eq.~(\ref{eq_multimodal}).}
\label{fig-HG4}
\end{center}
\end{figure}

In practice, phase functions are obtained from scattering calculations using Mie theory, see e.g. \cite{wiscombe1980}. Such phase functions may have multiple local extrema. An artificial example (not obtained from Mie theory) is given in  fig.~ \ref{fig-HG2b}. It uses the function 
\begin{eqnarray}
p(x) &=& 0.8 p_{HG}(x;.9) + 0.1p_{HG}(x;-.6)\\
 &\,& + 
0.04f_1(x; .2, .01 ,3) + 0.06f_2(x;.6 , .02) \, .
\label{eq_multimodal}
\end{eqnarray}
where $f_1$ and $f_2$ are as in eq.~(\ref{eq_f1}, \ref{eq_f2}). The integrand in eq.~(\ref{eq_matrix}) turns out to be analytic in the strip $\cS_\alpha$ with $\alpha \approx 10^{-2}$.  
The scattering kernels $P_0$ and $P_7$ for this phase function were computed at $200 \times 200$ points with eq.~(\ref{eq_trapezoidal}), using $N = 128$ terms in each case. A logarithmic heat map of $P_0$ is shown in fig.~\ref{fig-HG3} and a heat map of $P_7$ is shown in fig.~\ref{fig-HG4}. The calculation took about 10 seconds per scattering kernel on a laptop equipped with a dual-core processor running at 1.40 GHz. The relative accuracy of each result is about $10^{-3}$.     

\section{Conclusion}

A direct numerical integration method using the trapezoidal rule has been presented for the evaluation of scattering kernels that arise in plane-parallel radiation transfer equations. Its convergence is exponential, and the relation between the convergence rate and the domain of analyticity of the phase function is explained. The note also presents a closed form of the scattering kernel for the Henyey-Greenstein phase function, in terms of complete elliptic integrals of the second kind. The closed form can be used to assess the accuracy of the proposed numerical integration scheme. 

Most computational approaches to plane-parallel radiative transfer use discretizations based on truncated versions of eq.~(\ref{eq_LegendreExp}) (Nystr\"om's method). However, some problems of this form also require the evaluation of intensities  from scattered beams which may be computed from scattering kernels. This is where a fast and accurate computational scheme such as the one presented here will hopefully be of use.

\section{Appendix: Complete elliptic integrals} Legendre's complete elliptic integrals $K_0$ and $E_0$ of the first and second kind are defined as
\begin{eqnarray}
K_0(m) &=& \int_0^{\pi/2} \frac{ds}{\sqrt{1 - m \sin^2 s}} \\
&=& 
\frac{1}{2} \int_0^\pi \frac{ds}{\sqrt{1 - \frac{m}{2} \pm \frac{m}{2} \cos s}} 
\label{eq_K0}
\\ 
E_0(m) &=& \int_0^{\pi/2} \sqrt{1 - m \sin^2 s} \, ds \\
&=&
\frac{1}{2} \int_0^\pi \sqrt{1 - \frac{m}{2} \pm \frac{m}{2} \cos s} \, ds 
\, .
\label{eq_E0}
\end{eqnarray}   
where $m \in \C$ is known as the \emph{modulus}. Note that usually 
these integrals are expressed in terms of the \emph{parameter} $k$ where  $m = k^2$, leading to the more common notation
\begin{equation}
K(k) = K_0(k^2), \quad E(k) = E_0(k^2)
\label{eq-elliptic-parameter}
\end{equation}
For example, the treatment in \cite{nist2010} is in terms of $K, \, E$ while \textsc{Mathematica}\textsuperscript{\textregistered} uses $K_0, \, E_0$.
These integrals converge for $k \in \C$ with $\Re m < 1$ and the functions can be continued analytically to $\C$ minus a branch cut along $[1,\infty)$.  
 It is known (see eq.~(19.4.1) in \cite{nist2010}) that 
\begin{equation}
2m\frac{d}{dm} K_0(m) = \frac{E_0(m)}{1-m} - K_0(m) \, .
\label{eq-dK}
\end{equation}
Therefore, for $\alpha,\,  \beta \in \C$, we can set $m = \frac{2\beta}{\alpha + \beta}$ and obtain in the case when $\Re \frac{2\beta}{\alpha + \beta} < 1$
\begin{eqnarray}
\int_0^\pi \frac{ds}{\sqrt{\alpha - \beta \cos s}} &=& \frac{2}{\sqrt{\alpha + \beta}} K_0\left(\frac{2\beta}{\alpha + \beta} \right) 
\label{eq_K1}
\\
\int_0^\pi \sqrt{\alpha - \beta \cos s} \, ds &=& 2\sqrt{\alpha + \beta} E_0\left(
\frac{2\beta}{\alpha + \beta} \right) 
\label{eq_E1}
\end{eqnarray}  
where the principal branch of the square root is used.

Given real $0 < m < 1$, then $E_0(m)$ and $K_0(m)$ may be evaluated rapidly using the iterations
\begin{eqnarray}
a_0 &=& 1, \quad g_0 = \sqrt{1 - m}, \quad c_0 = \sqrt{m}\\
a_{n+1} &=& \frac{a_n + g_n}{2} \quad g_{n+1} = \sqrt{a_n g_n} \\
c_{n+1} &=& \frac{c_n^2}{4a_{n+1}} 
\end{eqnarray}
for $n = 0, \, 1, \, 2, \dots$.
Then 
\[\lim_{n \to \infty} a_n = \lim_{n \to \infty} g_n = M
\]
where $M = AGM(1,g_0)$ is known as Gauss's arithmetic-geometric mean; see \cite{nist2010}. The convergence is quadratic. Moreover, $\lim_{n \to \infty} c_n  = 0$, and the convergence is also quadratic. One then obtains the values of $K_0$ and $E_0$ from
\begin{equation}
K_0(m) = \frac{2}{\pi M}, \quad E_0(m) = \frac{2}{\pi M} \left(1 - \sum_{n = 0}^\infty 2^{n-1}c_n^2 \right)\, . 
\label{eq-fastelliptic}
\end{equation} 
Due to the rapid convergence, only a few terms need to be evaluated. An alternative fast computation method for $E_0(M)$ is given in \cite{adlaj2012}. 

\medskip
\textbf{Acknowledgement} This research was supported by the Cooperative Institute for Research in the Atmosphere (CIRA) at Colorado State University. Part of this work was carried out at the Joint Center for Satellite Data Assimilation (JCSDA) at NCWCP, College Park, MD.

\bibliographystyle{plain}  
\bibliography{hg}

\end{document}